\documentclass[12pt]{amsart}

\usepackage{amssymb,latexsym}

\usepackage{enumerate}

\usepackage[T1]{fontenc}

\numberwithin{equation}{section} 
\numberwithin{figure}{section} 
\theoremstyle{plain}
\theoremstyle{plain}
\newtheorem{thm}{Theorem}
  \theoremstyle{definition}
  \newtheorem{defn}[thm]{Definition}
 \theoremstyle{definition}
  \newtheorem{example}[thm]{Example}
  \theoremstyle{remark}
  \newtheorem{rem}[thm]{Remark}
  \theoremstyle{plain}
  \newtheorem{lem}[thm]{Lemma}
  \theoremstyle{plain}
  \newtheorem*{prop*}{Proposition}

\begin{document}

\baselineskip=17pt

\title[Stafney's lemma for ``classical" interpolation methods]{Stafney's lemma holds for several {}``classical'' interpolation
methods}
\author[A. Ivtsan]{Alon Ivtsan}

\address{Department of Mathematics, Technion I.I.T. Haifa 32000, Israel}

\email{aloniv@techunix.technion.ac.il}


\begin{abstract}
Let $\left(B_{0},B_{1}\right)$ be a Banach pair. Stafney showed that one can replace the space $\mathcal{F}$$\left(B_{0},B_{1}\right)$
by its dense subspace $\mathcal{G}\left(B_{0},B_{1}\right)$ 
in the definition of the norm in the Calderón complex interpolation
method on the strip if
the element belongs to the intersection of the spaces $B_{i}$. We
shall extend this result to a more general setting, which contains
well-known interpolation methods: the Calderón complex interpolation
method on the annulus, the Lions-Peetre real method (with several
different choices of norms), and the Peetre {}``$\pm$'' method.
\end{abstract}
\maketitle
\section{Introduction}
Stafney showed in his paper \cite{stafney} (Lemma 2.5, p. 335) that one can
replace the space $\mathcal{F}\left(B_{0},B_{1}\right)$ by its
dense subspace $\mathcal{G}\left(B_{0},B_{1}\right)$ 
in the definition of the norm of an element in the Calderón complex
interpolation space $\left[B_{0},B_{1}\right]_{\theta}$ if the element
belongs to the intersection of the two Banach spaces. Among the various
applications of Stafney's result, we mention that it can be used to
give an apparently simpler proof of part of Calderón's duality theorem,
namely that, for each $\theta\in\left(0,1\right)$, $\left[B_{0}^{*},B_{1}^{*}\right]^{\theta}\subset\left(\left[B_{0},B_{1}\right]_{\theta}\right)^{*}$
holds if $\left(B_{0},B_{1}\right)$ is a regular Banach couple (\cite{cwikel-lecture},
pp.~20-1). 

We shall obtain a version of Stafney's lemma in the general setting
of pseudolattices, which, for appropriate selections of the parameters,
will give us analogues of this lemma for the Calderón complex interpolation
method on the annulus, for the {}``discrete definition'' of the
Lions-Peetre real method and for the Peetre {}``$\pm$'' method.
In the case of the Lions-Peetre method we can work either with the
norm defined via the $J$-functional or with an earlier used variant
of that norm introduced in \cite{lions-peetre}. The formulation of
our version of Stafney's lemma for these methods can be found in Remark
\ref{rem:explicit-formulation}. 

There is some overlap of our result for the case of the Lions-Peetre
method with a result in Appendix 3 on pp.~47--48 of \cite{cwikel-peetre}.
The latter result applies to more general versions of interpolation
spaces defined by the $K$-functional, but, on the other hand, the
relevant norm estimate there is only to within equivalence of norms.
See also the {}``Note added in proof'' on p.~49 of \cite{cwikel-peetre},
which announces (without explicit proof) that the condition of {}``mutual
closedness'' imposed in Appendix 3 can be removed.

Stafney's lemma cannot be extended to all interpolation methods where
there are natural analogues of the space $\mathcal{F}$ and its dense
subspace $\mathcal{G}$. For example, in \cite{ccrsw} an analogue
of Calderón's complex interpolation space generated by an $n$-tuple
of Banach spaces, rather than just a couple, is introduced. The example
in Appendix 1 on pp.~223--6 of \cite{ccrsw} shows that, perhaps
surprisingly, for $n\ge3$, the expected analogue of Stafney's lemma
can fail to hold. In fact, in that setting the natural analogues of
the two quantities which appear below in the formula (\ref{eq:mr})
in Theorem \ref{thm:main} may even fail to be equivalent.

\section{Preliminaries and examples}

Before stating our main result we need to provide a number of definitions
and examples, most of which are from \cite{ckmr}.
\begin{defn}
Let $Ban$ be the class of all Banach spaces over the complex numbers.
A mapping $\mathcal{X}:Ban\rightarrow Ban$ will be called a \emph{pseudolattice}
if

(i) for each $B\in Ban$ the space $\mathcal{X}\left(B\right)$ consists
of $B$ valued sequences $\left\{ b_{n}\right\} _{n\in\mathbb{Z}}$, 

(ii) whenever $A$ is a closed subspace of $B$ it follows that $\mathcal{X}\left(A\right)$
is a closed subspace of $\mathcal{X}\left(B\right)$, and 

(iii) there exists a positive constant $C=C\left(\mathcal{X}\right)$
such that, for all $A,B\in Ban$ and all bounded linear operators
$T:A\rightarrow B$ and every sequence $\left\{ a_{n}\right\} _{n\in\mathbb{Z}}\in\mathcal{X}\left(A\right)$,
the sequence $\left\{ Ta_{n}\right\} _{n\in\mathbb{Z}}\in\mathcal{X}\left(B\right)$
and satisfies the estimate \[
\left\Vert \left\{ Ta_{n}\right\} _{n\in\mathbb{Z}}\right\Vert _{\mathcal{X}\left(B\right)}\leq C\left(\mathcal{X}\right)\left\Vert T\right\Vert _{A\rightarrow B}\left\Vert \left\{ a_{n}\right\} _{n\in\mathbb{Z}}\right\Vert _{\mathcal{X}\left(A\right)}\,.\]

\end{defn}
The following examples will be relevant for our applications.
\begin{example}
Let $X$ be a Banach lattice of real valued functions defined on $\mathbb{Z}$.
We will use the notation $\mathcal{X}=X$ to mean that, for each $B\in Ban$,
$\mathcal{X}\left(B\right)$ is the space, usually denoted by $X\left(B\right)$,
consisting of all $B$ valued sequences $\left\{ b_{n}\right\} _{n\in\mathbb{Z}}$
such that $\left\{ \left\Vert b_{n}\right\Vert _{B}\right\} _{n\in\mathbb{Z}}\in X$.
It is normed by $\left\Vert \left\{ b_{n}\right\} _{n\in\mathbb{Z}}\right\Vert _{X\left(B\right)}=\left\Vert \left\{ \left\Vert b_{n}\right\Vert _{B}\right\} _{n\in\mathbb{Z}}\right\Vert _{X}$.
In particular, we shall be interested in the choices $X=\ell^{p}$
for $p\in[1,\infty]$ and $X=c_{0}$.
\end{example}

\begin{example}
For each $B\in Ban$ let $FC\left(B\right)$ be the space of all $B$
valued sequences $\left\{ b_{n}\right\} _{n\in\mathbb{Z}}$ such that
$b_{n}=\frac{1}{2\pi}\int_{0}^{2\pi}e^{-int}f\left(e^{it}\right)dt$
for all $n$ and some continuous function $f:\mathbb{T}\rightarrow B$.
$FC\left(B\right)$ is normed by $\left\Vert \left\{ b_{n}\right\} _{n\in\mathbb{Z}}\right\Vert _{FC\left(B\right)}=\sup_{t\in[0,2\pi)}\left\Vert f\left(e^{it}\right)\right\Vert _{B}$.
The notation $\mathcal{X}=FC$ will mean that $\mathcal{X}\left(B\right)=FC\left(B\right)$
for each $B$.
\end{example}

\begin{example}
\label{exa:pm}We shall use the notation $\mathcal{X}=UC$, when $\mathcal{X}\left(B\right)=UC\left(B\right)$
for every $B\in Ban$, where $UC\left(B\right)$ denotes the Banach
space of all $B$ valued sequences $\left\{ b_{n}\right\} _{n\in\mathbb{Z}}$
such that $\sum_{n\in\mathbb{Z}}\lambda_{n}b_{n}$ converges for all
complex sequences $\left\{ \lambda_{n}\right\} _{n\in\mathbb{Z}}$
satisfying $\left|\lambda_{n}\right|\leq1$ for each $n\in\mathbb{Z}$
(i.e. such that the sequence $b_{n}$ is unconditionally convergent),
and $UC(B)$ is normed by $\left\Vert \left\{ b_{n}\right\} _{n\in\mathbb{Z}}\right\Vert _{UC\left(B\right)}=\sup\left\{ \left\Vert \sum_{n\in F}\lambda_{n}b_{n}\right\Vert _{B}\right\} $,
where the supremum is taken over all finite subsets $F$ of $\mathbb{Z}$
and all sequences $\left\{ \lambda_{n}\right\} _{n\in\mathbb{Z}}$
which satisfy $\left|\lambda_{n}\right|\leq1$ for all $n$ (see \cite{peetre}
pp.~174-5). (Note that, as was pointed out on p.~58 of \cite{janson},
it suffices to consider sequences $\left\{ \lambda_{n}\right\} _{n\in\mathbb{Z}}$
with $\lambda_{n}=\pm1$, since this yields the same space to within
equivalence of norms.)

Analogously, we shall use the notation $\mathcal{X}=WUC$, when $\mathcal{X}\left(B\right)=WUC\left(B\right)$
for every $B\in Ban$, where $WUC\left(B\right)$ denotes the space
consisting of all $B$ valued sequences $\left\{ b_{n}\right\} _{n\in\mathbb{Z}}$
for which the above norm $\left\Vert \left\{ b_{n}\right\} _{n\in\mathbb{\mathbb{Z}}}\right\Vert _{UC(B)}$
is finite, but for which the unconditional convergence of the sequence
$b_{n}$ is not required. Such sequences are said to be weakly unconditionally
convergent. (See e.g. p.~58 of \cite{janson} or pp.~99 of \cite{lind-tzaf}
for an equivalent definition.) 
\end{example}

Let $\overline{B}=\left(B_{0},B_{1}\right)$ be a Banach pair (i.e.
$B_{0}$ and $B_{1}$ are two Banach spaces which are continuously
embedded in some Hausdorff topological vector space). Let $\mathcal{X}_{0}$
and $\mathcal{X}_{1}$ be any two pseudolattices. We consider them
as a pair, which we denote by $\mathbf{X}=\left\{ \mathcal{X}_{0},\mathcal{X}_{1}\right\} $. 
\begin{defn}
\label{def:j_x_b}For each Banach pair $\bar{B}$ and pseudolattice
pair $\mathbf{X}$ we define $\mathcal{J}\left(\mathbf{X},\bar{B}\right)$
to be the space of all $B_{0}\cap B_{1}$ valued sequences $\left\{ b_{n}\right\} _{n\in\mathbb{Z}}$
for which the sequence $\left\{ e^{jn}b_{n}\right\} _{n\in\mathbb{Z}}$
is in $\mathcal{X}_{j}\left(B_{j}\right)$ for $j=0,1$. This space
is normed by \[
\left\Vert \left\{ b_{n}\right\} _{n\in\mathbb{Z}}\right\Vert _{\mathcal{J}\left(\mathbf{X},\bar{B}\right)}=\max_{j=0,1}\left\Vert \left\{ e^{jn}b_{n}\right\} _{n\in\mathbb{Z}}\right\Vert _{\mathcal{X}_{j}\left(B_{j}\right)}\,.\]

\end{defn}

\begin{defn}
\label{def:nontrivial}Let $\mathbb{A}$ denote the annulus $\left\{ z\in\mathbb{C}:1\leq\left|z\right|\leq e\right\} $
and let $\mathbb{A}^{\circ}$ denote its interior. We shall say that
the pseudolattice pair $\mathbf{X}$ is \emph{nontrivial}
if, for the special one-dimensional Banach pair $\bar{B}=\left(\mathbb{C},\mathbb{C}\right)$
and each $s\in\mathbb{A}^{\circ}$, there exists a sequence $\left\{ b_{n}\right\} _{n\in\mathbb{Z}}\in\mathcal{J}\left(\mathbf{X},\bar{B}\right)$
such that the limit $\lim{}_{M,N\to+\infty}\sum_{n=-M}^{N}s^{n}b_{n}$
exists and is finite and non zero.
\end{defn}
(Note that here our notation differs slightly from that in \cite{ckmr},
but is the same as that in \cite{cwikel-janson}.)
\begin{defn}
\label{def:laurent_compatible}We shall say that the pseudolattice pair $\mathbf{X}$ is \emph{Laurent compatible} if it is
nontrivial and if for every Banach pair $\bar{B}$, every sequence
$\left\{ b_{n}\right\} _{n\in\mathbb{Z}}$ in $\mathcal{J}\left(\mathbf{X},\bar{B}\right)$
and every fixed $z$ in the open annulus $\mathbb{A}^{\circ}$, the
Laurent series $\sum_{n\in\mathbb{Z}}z^{n}b_{n}$ converges in $B_{0}+B_{1}$
and $\left\Vert \sum_{n\in\mathbb{Z}}z^{n}b_{n}\right\Vert _{B_{0}+B_{1}}\leq C\left\Vert \left\{ b_{n}\right\} _{n\in\mathbb{Z}}\right\Vert _{\mathcal{J}\left(\mathbf{X},\bar{B}\right)}$
for some constant $C=C\left(z\right)$ independent of the choice of
$\left\{ b_{n}\right\} _{n\in\mathbb{Z}}$.\end{defn}
\begin{rem}
\label{rem:analytic}If $\mathbf{X}$ is Laurent compatible, then
the $B_{0}+B_{1}$ valued function $f\left(z\right)=\sum_{n\in\mathbb{Z}}z^{n}b_{n}$
is analytic in $\mathbb{A}^{\circ}$. \end{rem}
\begin{defn}
\label{def:b_x_s}For each Banach pair $\bar{B}$, each Laurent compatible
pair $\mathbf{X}$ and each fixed $s\in\mathbb{A}^{\circ}$ we define
the space $\bar{B}_{\mathbf{X},s}$ to consist of all the elements
of the form $b=\sum_{n\in\mathbb{Z}}s^{n}b_{n}$ where $\left\{ b_{n}\right\} _{n\in\mathbb{Z}}\in\mathcal{J}\left(\mathbf{X},\bar{B}\right)$,
with the norm $\left\Vert b\right\Vert _{\bar{B}_{\mathbf{X},s}}=\inf\left\{ \left\Vert \left\{ b_{n}\right\} _{n\in\mathbb{Z}}\right\Vert _{\mathcal{J}\left(\mathbf{X},\bar{B}\right)}:b=\sum_{n\in\mathbb{Z}}s^{n}b_{n}\right\} $. \end{defn}
\begin{rem}
$\bar{B}_{\mathbf{X},s}$ equipped with this norm is a Banach space.
\end{rem}

\begin{rem}
\label{rem:concrete examples}As remarked on p.~251 of \cite{ckmr},
in elaboration of a point of view going back to \cite{peetre}, the
space $\bar{B}_{\mathbf{X},s}$ coincides with various known interpolation
spaces for appropriate choices of $\mathcal{X}_{0}$, $\mathcal{X}_{1}$
and $s$. In each of the three following examples we set $s=e^{\theta}$
for some $\theta\in\left(0,1\right)$. 
\end{rem}
(i) If $\mathcal{X}_{0}=\mathcal{X}_{1}=FC$, the space $\bar{B}_{\mathbf{X},s}$
coincides isometrically with the variant of Calderón's complex interpolation
space obtained when an annulus is used instead of a strip in the interpolation
method. (The proof of this last claim is rather straightforward, but,
since it is not given explicitly in \cite{peetre} or \cite{ckmr},
we provide it for the reader's convenience in Section \ref{sec:Additional-remarks-regarding},
together with the relevant definitions.) Note that this variant of
Calderón's complex interpolation space, which was apparently first
considered in \cite{peetre}, coincides with $\left[B_{0},B_{1}\right]_{\theta}$
to within equivalence of norms, as was shown on pp.~1007-9 of \cite{cwikel-complex}.

(ii) If $\mathcal{X}_{0}=\mathcal{X}_{1}=\ell^{p}$, then $\bar{B}_{\mathbf{X},s}$
is the Lions-Peetre real method space $\bar{B}_{\theta,p}=\left(B_{0},B_{1}\right)_{\theta,p}$.
In fact the norm that we obtain here is exactly the norm introduced
in formula (1.3) on p.~17 of \cite{lions-peetre} for suitable choices
of the parameters $p_{0},$ $p_{1}$, $\xi_{0}$ and $\xi_{1}$. In
\cite{lions-peetre} this space is denoted by $s(p_{0},\xi_{0},B_{0};p_{1},\xi_{1},B_{1})$, a notation which is now rarely used. It has become more customary
to use other equivalent norms on $\left(B_{0},B_{1}\right)_{\theta,p}$
which are defined via the Peetre $J$-functional or $K$-functional.
For example, in terms of the $J$-functional (i.e. $J\left(t,x;B_{0},B_{1}\right)=\max_{j=0,1}t^{j}\left\Vert x\right\Vert _{B_{j}}$
for $x\in B_{0}\cap B_{1}$), one can use the norm
{$\left\Vert x\right\Vert =\inf\left\Vert \left\{ e^{-\theta n}J(e^{n},c_{n};B_{0},B_{1})\right\} _{n\in\mathbb{Z}}\right\Vert _{\ell^{p}}$,
where the infimum is taken over all representations $x=\sum_{n\in\mathbb{Z}}c_{n}$
(with convergence in $B_{0}+B_{1}$) with $\left\{ c_{n}\right\} _{n\in\mathbb{Z}}\in\mathcal{J}\left(\left\{ \ell^{p},\ell^{p}\right\} ,\bar{B}\right)$.}

(iii) If $\mathcal{X}_{0}=\mathcal{X}_{1}=UC$, then $\bar{B}_{\mathbf{X},s}$
is an appropriate slight modification (using powers of $e$ instead
of powers of $2$) of the interpolation space $\bar{B}_{\left\langle \theta\right\rangle }=\left\langle B_{0},B_{1}\right\rangle _{\theta}$
introduced by Peetre on p.~175--6 of \cite{peetre}, for the function
parameter $\rho(t)=t^{\theta}$, and if $\mathcal{X}_{0}=\mathcal{X}_{1}=WUC$,
then $\bar{B}_{\mathbf{X},s}$ is (a modification, again with powers
of $e$ in place of $2$ of) the Gustavsson-Peetre variant of $\left\langle B_{0},B_{1}\right\rangle _{\theta}$
which is denoted by $\left\langle \bar{B},\rho_{\theta}\right\rangle $
(see p.~45 of \cite{gust-peetre}). In fact, for these two choices of
the pseudolattice pair, the spaces $\bar{B}_{\mathbf{X},s}$ will
coincide exactly with the spaces $\left\langle B_{0},B_{1}\right\rangle _{\theta}$
and $\left\langle \bar{B},\rho_{\theta}\right\rangle $ respectively
if we perform the following rather obvious changes in their construction:
replace powers of $e$ by powers of $2$ in Definition \ref{def:j_x_b},
replace the annulus $\mathbb{A}$ by the annulus $\left\{ z\in\mathbb{C}:1\leq\left|z\right|\leq2\right\} $
in Definition \ref{def:b_x_s} and set $s=2^{\theta}$ for some $\theta\in\left(0,1\right)$.

(The method which yields the spaces $\left\langle B_{0},B_{1}\right\rangle _{\theta}$
is sometimes referred to as the {}``$\pm$'' method, since, as mentioned
above, in the definition of unconditional convergence it suffices
to consider sequences $\lambda_{n}$ whose values are $1$ and $-1$.)

\begin{defn}
Let $\mathbf{X}=\left\{ \mathcal{X}_{0},\mathcal{X}_{1}\right\} $
be a pair of pseudolattices. We shall say that $\mathbf{X}$ \emph{admits
differentiation} if it is Laurent compatible and, for each complex
Banach space $B$,

(i) for each $r\in\left(0,1\right)$, each element $\left\{ b_{n}\right\} _{n\in\mathbb{Z}}\in\mathcal{X}_{0}\left(B\right)$ satisfies 

$\lim_{k\rightarrow-\infty}r^{-k}\left\Vert b_{k}\right\Vert _{B}=0$
and each element $\left\{ b_{n}\right\} _{n\in\mathbb{Z}}\in\mathcal{X}_{1}\left(B\right)$ satisfies $\lim_{k\rightarrow\infty}r^{k}\left\Vert b_{k}\right\Vert _{B}=0$, and

(ii) for every complex number $\rho$ satisfying $0<\left|\rho\right|<1$ and
for every sequence $\left\{ b_{n}\right\} _{n\in\mathbb{Z}}\in\mathcal{X}_{0}\left(B\right)\cap\mathcal{X}_{1}\left(B\right)$,
the new sequence $\left\{ b_{n}^{j}\right\} _{n\in\mathbb{Z}}$ is
also in $\mathcal{X}_{j}\left(B\right)$ for $j=0,1$, where $\left\{ b_{n}^{0}\right\} _{n\in\mathbb{Z}}$
and $\left\{ b_{n}^{1}\right\} _{n\in\mathbb{Z}}$ are defined by
setting $b_{n}^{0}=\sum_{k<0}\rho^{-k}b_{n+k+1}$ and $b_{n}^{1}=\sum_{k\geq0}\rho^{k}b_{n+k+1}$
(where the convergence of these sums in $B$ is guaranteed by condition
(i)), and if also

(iii) for $j=0,1$ and each $\rho$ as above, the linear map $D_{j,\rho}$
defined on $\mathcal{X}_{j}\left(B\right)$ by setting $D_{j,\rho}\left(\left\{ b_{n}\right\} _{n\in\mathbb{Z}}\right)=\left\{ b_{n}^{j}\right\} _{n\in\mathbb{Z}}$
maps $\mathcal{X}_{j}\left(B\right)$ boundedly into itself.\end{defn}
\begin{rem}
\label{rem:admits_differentiation}The pair $\mathbf{X}=\left\{ \mathcal{X}_{0},\mathcal{X}_{1}\right\} $
admits differentiation whenever $\mathcal{X}_{0}$ and $\mathcal{X}_{1}$
are each chosen to be any of $\ell^{p}$ $(p\in[1,\infty])$, $c_{0}$,
$FC$, $UC$ or $WUC$ (see p. 256 of \cite{ckmr}).
\end{rem}
The property of admitting differentiation has the following consequence
(which also explains the choice of terminology for this property).
\begin{lem}
\label{lem:Lemma}Let $\mathbf{X}$ be a pair of pseudolattices which
admits differentiation and let $\bar{B}$ be a Banach pair. Let the
sequence $\left\{ f_{n}\right\} _{n\in\mathbb{Z}}$ be an element
of $\mathcal{J}\left(\mathbf{X},\bar{B}\right)$ and let $f:\mathbb{A}^{\circ}\rightarrow B_{0}+B_{1}$
be the analytic function defined by $f\left(z\right)=\sum_{n\in\mathbb{Z}}z^{n}f_{n}$.
Suppose that $f\left(s\right)=0$ for some point $s\in\mathbb{A}^{\circ}$
and let $g:\mathbb{A}^{\circ}\rightarrow B_{0}+B_{1}$ be the analytic
function obtained by setting $g\left(s\right)=f'\left(s\right)$ and
$g\left(z\right)=\frac{1}{z-s}f\left(z\right)$ for all $z\in\mathbb{A}^{\circ}\backslash\left\{ s\right\} $.
Let $\left\{ g_{n}\right\} _{n\in\mathbb{Z}}$ be the sequence of
coefficients in the Laurent expansion $g\left(z\right)=\sum_{n\in\mathbb{Z}}z^{n}g_{n}$
of $g$ in $\mathbb{A^{\circ}}$. Then $\left\{ g_{n}\right\} _{n\in\mathbb{Z}}$
is also an element of $\mathcal{J}\left(\mathbf{X},\bar{B}\right)$.
\end{lem}
For the proof we refer the reader to pp.~258-9 of \cite{ckmr}. This
lemma will also hold if we replace $e$ by some $r>1$ in the definitions
of $\mathbb{A}$ and $\mathcal{J}\left(\mathbf{X},\bar{B}\right)$.

We conclude this section with two more definitions of notions which
will appear explicitly in our main theorem.
\begin{defn}
For each Banach pair $\bar{B}$ we define $\mathcal{J}_{0}\left(\bar{B}\right)$
to be the space of all $B_{0}\cap B_{1}$ valued sequences $\left\{ b_{n}\right\} _{n\in\mathbb{Z}}$
with finite support.\end{defn}
\begin{rem}
\label{rem:j_zero}For every Banach pair $\bar{B}$, we obviously
have that \[
\mathcal{J}_{0}\left(\bar{B}\right)\subset\mathcal{J}\left(\left\{ \mathcal{X}_{0},\mathcal{X}_{1}\right\} ,\bar{B}\right)\]
 whenever $\mathcal{X}_{0}$ and $\mathcal{X}_{1}$ are chosen to
be any of the pseudolattices $FC$, $UC$, $WUC$, $\ell^{p}$ for
$p\in[1,\infty]$ or $c_{0}$. Furthermore, one can verify that $\mathcal{J}_{0}\left(\bar{B}\right)$
is dense in $\mathcal{J}\left(\left\{ FC,FC\right\} ,\bar{B}\right)$
and also dense in $\mathcal{J}\left(\left\{ \mathcal{X}_{0},\mathcal{X}_{1}\right\} ,\bar{B}\right)$
for $\mathcal{X}_{i}\in\left\{ UC,c_{0},\ell^{p}\right\} $, $1\leq p<\infty$,
$i=0,1$. But, in general, $\mathcal{J}_{0}\left(\bar{B}\right)$
is not dense in $\mathcal{J}\left(\left\{ WUC,WUC\right\} ,\bar{B}\right)$
and, except for trivial Banach spaces $B_{0},B_{1}$, is never dense
in $\mathcal{J}\left(\left\{ \ell^{\infty},\ell^{\infty}\right\} ,\bar{B}\right)$. 
\end{rem}

\begin{defn}
Let $S$ denote the right-shift operator on two-sided sequences defined
by $S\left(\left\{ b_{n}\right\} _{n\in\mathbb{Z}}\right)=\left\{ b_{n-1}\right\} _{n\in\mathbb{Z}}$
. 
\end{defn}

\section{The main theorem}

We can now state and prove our main theorem. The steps of the proof
parallel the steps of Stafney's proof on p.~335 of \cite{stafney}.
\begin{thm}
\label{thm:main}Let $\mathbf{X}$ be a pair of pseudolattices which
admits differentiation and let $\bar{B}$ be a Banach pair. Suppose
that 

(i) $\mathcal{J}_{0}\left(\bar{B}\right)\subset\mathcal{J}\left(\mathbf{X},\bar{B}\right)$
and $\mathcal{J}_{0}\left(\bar{B}\right)$ is dense in $\mathcal{J}\left(\mathbf{X},\bar{B}\right)$
and that 

(ii) The right-shift operator $S$ maps $\mathcal{X}_{j}\left(B_{j}\right)$
boundedly into itself for $j=0,1$. 

Then, for each $x\in B_{0}\cap B_{1}$ and $s\in\mathbb{A}^{\circ}$, 

\begin{equation}
\left\Vert x\right\Vert _{\bar{B}_{\mathbf{X},s}}=\inf\left\{ \left\Vert \left\{ b_{n}\right\} _{n\in\mathbb{Z}}\right\Vert _{\mathcal{J}\left(\mathbf{X},\bar{B}\right)}:\sum_{n\in\mathbb{Z}}s^{n}b_{n}=x,\ \left\{ b_{n}\right\} _{n\in\mathbb{Z}}\in\mathcal{J}_{0}\left(\bar{B}\right)\right\} \,.\label{eq:mr}\end{equation}

\end{thm}

\begin{rem}
\label{rem:equivalent}Our proof will also hold if we replace the
norm $\left\Vert \cdot\right\Vert _{\mathcal{J}\left(\mathbf{X},\bar{B}\right)}$
in both Definition \ref{def:b_x_s} and Equation (\ref{eq:mr}) by
an equivalent one.
\end{rem}

\begin{rem}
\label{rem:replace}As hinted in part (iii) of Remark \ref{rem:concrete examples}
and in the remark following the statement of Lemma \ref{lem:Lemma},
if we replace $e$ in our definitions by any positive number greater
than $1$, then we can obtain an appropriate reformulation of Theorem
\ref{thm:main}.
\end{rem}

\begin{rem}
\label{rem:explicit-formulation}By Remarks \ref{rem:concrete examples}, \ref{rem:admits_differentiation}, and \ref{rem:j_zero} we can obtain an
appropriate formulation of Theorem \ref{thm:main} for the {}``annulus''
variant of Calderón's complex interpolation method space, for the
Lions-Peetre real method space $\left(B_{0},B_{1}\right)_{\theta,p}$
for $1\leq p<\infty$ and for the Peetre interpolation space $\bar{B}_{\left\langle \theta\right\rangle }=\left\langle B_{0},B_{1}\right\rangle _{\theta}$
for $\theta\in\left(0,1\right)$. (Of course condition (ii) of Theorem
\ref{thm:main} obviously holds in these cases and in fact $S$ is
even an isometry.)

Our theorem in the case of the {}``annulus'' variant of Calderón's
complex interpolation method space can also be obtained by an alternative
argument similar to Stafney's proof on p.~335 of \cite{stafney},
if one replaces Calderón's space $\mathcal{G}\left(B_{0},B_{1}\right)$
by the space of all Laurent polynomials with coefficients in $B_{0}\cap B_{1}$
and if $\mathcal{F}\left(B_{0},B_{1}\right)$ is replaced by its counterpart
for the annulus (see also Section \ref{sec:Additional-remarks-regarding}
and the paragraph which precedes Definition 4.1 on p.~80 of \cite{cwikel-janson}). 

For the Lions-Peetre real method space $\left(B_{0},B_{1}\right)_{\theta,p}$
for $1\leq p<\infty$, our theorem shows that if $x\in B_{0}\cap B_{1}$,
$\theta\in\left(0,1\right)$ and $p\in[1,\infty)$, then 
\begin{eqnarray*}
 &  & \left\Vert x\right\Vert _{\left(B_{0},B_{1}\right)_{\theta,p}}\\
 & = & \inf\left\{ \max_{j=0,1}\left\Vert \left\{ e^{\left(j-\theta\right)n}c_{n}\right\} _{n\in\mathbb{Z}}\right\Vert _{\ell^{p}\left(B_{j}\right)}:x=\sum_{n\in\mathbb{Z}}c_{n},\left\{ c_{n}\right\} _{n\in\mathbb{Z}}\in\mathcal{J}_{0}\left(\bar{B}\right)\right\} \,.\end{eqnarray*}
By Remark \ref{rem:equivalent}, if we equip the space $\left(B_{0},B_{1}\right)_{\theta,p}$
with the norm {$\left\Vert x\right\Vert =\inf\left\Vert \left\{ e^{-\theta n}J(e^{n},c_{n};B_{0},B_{1})\right\} _{n\in\mathbb{Z}}\right\Vert _{\ell^{p}}$,
where the infimum is taken over all representations $x=\sum_{n\in\mathbb{Z}}c_{n}$
(with convergence in $B_{0}+B_{1}$) with $\left\{ c_{n}\right\} _{n\in\mathbb{Z}}\in\mathcal{J}\left(\left\{ \ell^{p},\ell^{p}\right\} ,\bar{B}\right)$}, then our theorem shows that for every $x\in B_{0}\cap B_{1}$, $\theta\in\left(0,1\right)$
and $p\in[1,\infty)$, we have that \textbf{\[
\left\Vert x\right\Vert =\inf\left\{ \left\Vert \left\{ e^{-\theta n}J(e^{n},c_{n};B_{0},B_{1})\right\} _{n\in\mathbb{Z}}\right\Vert _{\ell^{p}}:x=\sum_{n\in\mathbb{Z}}c_{n},\left\{ c_{n}\right\} _{n\in\mathbb{Z}}\in\mathcal{J}_{0}\left(\bar{B}\right)\right\}.\]}

By Remark \ref{rem:replace}, we can also obtain a version of our
theorem if, for instance, we equip the space $\left(B_{0},B_{1}\right)_{\theta,p}$ with the $J$-functional norm $\left\Vert x\right\Vert =\inf\left\Vert \left\{ 2^{-\theta n}J(2^{n},c_{n};B_{0},B_{1})\right\} _{n\in\mathbb{Z}}\right\Vert _{\ell^{p}}$, where the infimum is taken over all representations $x=\sum_{n\in\mathbb{Z}}c_{n}$ (with convergence in $B_{0}+B_{1}$) with $\left\{ c_{n}\right\}_{n\in\mathbb{Z}}$ belonging to the variant of $\mathcal{J}\left(\left\{ \ell^{p},\ell^{p}\right\} ,\bar{B}\right)$ obtained by replacing powers of $e$ by powers of $2$ (this norm appears on p.~43 of \cite{bergh-lofstrom}).

If we choose this norm, then for every $x\in B_{0}\cap B_{1}$, $\theta\in\left(0,1\right)$
and $p\in[1,\infty)$ we obtain that \textbf{\[
\left\Vert x\right\Vert =\inf\left\{ \left\Vert \left\{ 2^{-\theta n}J(2^{n},c_{n};B_{0},B_{1})\right\} _{n\in\mathbb{Z}}\right\Vert _{\ell^{p}}:x=\sum_{n\in\mathbb{Z}}c_{n},\left\{ c_{n}\right\} _{n\in\mathbb{Z}}\in\mathcal{J}_{0}\left(\bar{B}\right)\right\}.\]
}

By part (iii) of Remark \ref{rem:concrete examples} and by Remark
\ref{rem:replace}, for the Peetre interpolation space $\bar{B}_{\left\langle \theta\right\rangle }=\left\langle B_{0},B_{1}\right\rangle _{\theta}$
for $\theta\in\left(0,1\right)$, our theorem shows that if $x\in B_{0}\cap B_{1}$,
then
\begin{eqnarray*}
 &  & \left\Vert x\right\Vert _{\left\langle B_{0},B_{1}\right\rangle _{\theta}}\\
 & = & \inf\left\{ \max_{j=0,1}\left\Vert \left\{ 2^{\left(j-\theta\right)n}c_{n}\right\} _{n\in\mathbb{Z}}\right\Vert _{UC\left(B_{j}\right)}:x=\sum_{n\in\mathbb{Z}}c_{n},\left\{ c_{n}\right\} _{n\in\mathbb{Z}}\in\mathcal{J}_{0}\left(\bar{B}\right)\right\} \,.\end{eqnarray*}\end{rem}
\begin{rem}
We can in fact also obtain a version of our theorem for a discrete
version of the (generalised) $\mathcal{J}$-method which is discussed
(for example) on p.~381 of \cite{brudnyi} and apparently originated
in the work of Peetre in \cite{peetre-brazil}. We shall recall its
definition.
\end{rem}
Let $\Phi$ be a Banach lattice of two sided sequences satisfying
$\left\{ 0\right\} \neq\Phi\subset\ell_{1}^{0}+\ell_{1}^{1}$. Here
$\ell_{1}^{j}$ denotes the space of all real valued sequences $\left\{ b_{n}\right\} _{n\in\mathbb{Z}}$
such that the sum $\sum_{n\in\mathbb{Z}}2^{-nj}\left|b_{n}\right|$
is finite with the norm $\left\Vert \left\{ b_{n}\right\} _{n\in\mathbb{Z}}\right\Vert _{\ell_{1}^{j}}=\sum_{n\in\mathbb{Z}}2^{-nj}\left|b_{n}\right|$
for $j=0,1$. 

We define $\left\Vert x\right\Vert _{J_{\Phi}^{d}\left(\bar{B}\right)}=\inf\left\Vert \left\{ J\left(2^{n},x_{n};B_{0},B_{1}\right)\right\} _{n\in\mathbb{Z}}\right\Vert _{\Phi}$
where the infimum is taken over all representations $x=\sum_{n\in\mathbb{Z}}x_{n}$
(with convergence in $B_{0}+B_{1}$) with $x_{n}\in B_{0}\cap B_{1}$. 

We shall now show what is required in order to obtain our theorem
for the discrete $\mathcal{J}$-method. Set $x_{n}=2^{\theta n}c_{n}$.
Defining $\left\Vert \left\{ b_{n}\right\} _{n\in\mathbb{Z}}\right\Vert _{X}=\left\Vert \left\{ 2^{\theta n}b_{n}\right\} _{n\in\mathbb{Z}}\right\Vert _{\Phi}$
yields $\left\Vert x\right\Vert _{J_{\Phi}^{d}\left(\bar{B}\right)}=\inf\left\Vert \left\{ J\left(2^{n},c_{n};B_{0},B_{1}\right)\right\} _{n\in\mathbb{Z}}\right\Vert _{X}$.

Since $\left\Vert \left\{ J\left(2^{n},c_{n};B_{0},B_{1}\right)\right\} _{n\in\mathbb{Z}}\right\Vert _{X}$
and \[
\max_{j=0,1}\left\Vert \left\{ 2^{nj}\left\Vert c_{n}\right\Vert _{B_{j}}\right\} _{n\in\mathbb{Z}}\right\Vert _{X}=\max_{j=0,1}\left\Vert \left\{ 2^{nj}c_{n}\right\} _{n\in\mathbb{Z}}\right\Vert _{X\left(B_{j}\right)}\]
are equivalent, by Remarks \ref{rem:equivalent} and \ref{rem:replace}
we can obtain a version of Stafney's lemma in this case if the pseudolattice
pair $\left\{ \mathcal{X}_{0},\mathcal{X}_{1}\right\} =\left\{ X,X\right\} $
satisfies the conditions of our main theorem. (Here, of course, we
must replace $e$ by $2$ in the appropriate definitions.)

\begin{rem}
Janson showed that if we equip $B_{0}\cap B_{1}$ with the norms of
$\left\langle \bar{B},\rho_{\theta}\right\rangle $ and $\bar{B}_{\left\langle \theta\right\rangle }$,
we obtain two normed spaces with equivalent norms (see pp.~59-60 of\cite{janson}),
and thus a weaker version of (\ref{eq:mr}), i.e. that
the left and right sides are equivalent, can also be obtained for
$\mathcal{X}_{0}=\mathcal{X}_{1}=WUC$, even though, as pointed out
in Remark \ref{rem:j_zero}, condition (i) fails to hold in this case.
\end{rem}
\textit{Proof of the theorem.} Let $x$ be in $B_{0}\cap B_{1}$,
$s$ in $\mathbb{A}^{\circ}$ and $\varepsilon$ an arbitrary positive
number. The sequence $\left\{ b_{n}\right\} _{n\in\mathbb{Z}}$ defined
by setting $b_{0}=x$ and $b_{n}=0$ for $n\neq0$ is in $\mathcal{J}_{0}\left(\bar{B}\right)$
and satisfies $\sum_{n\in\mathbb{Z}}s^{n}b_{n}=x$. It is clear from
the definition of the norm $\left\Vert \cdot\right\Vert _{\bar{B}_{\mathbf{X},s}}$
that, for some $\left\{ c_{n}\right\} _{n\in\mathbb{Z}}$ in the subspace
$\mathcal{N}_{s}\left(\mathbf{X},\bar{B}\right)$ of $\mathcal{J}\left(\mathbf{X},\bar{B}\right)$,
consisting of all sequences $\left\{ d_{n}\right\} _{n\in\mathbb{Z}}$
such that $\sum_{n\in\mathbb{Z}}s^{n}d_{n}=0$, $\left\Vert \left\{ b_{n}\right\} _{n\in\mathbb{Z}}-\left\{ c_{n}\right\} _{n\in\mathbb{Z}}\right\Vert _{\mathcal{J}\left(\mathbf{X},\bar{B}\right)}<\left\Vert x\right\Vert _{\bar{B}_{\mathbf{X},s}}+\varepsilon/2$.
We need the following proposition:
\begin{prop*}
$\mathcal{J}_{0}\left(\bar{B}\right)\cap\mathcal{N}_{s}\left(\mathbf{X},\bar{B}\right)$
is dense in $\mathcal{N}_{s}\left(\mathbf{X},\bar{B}\right)$ with
respect to the norm of $\mathcal{J}\left(\mathbf{X},\bar{B}\right)$
restricted to $\mathcal{N}_{s}\left(\mathbf{X},\bar{B}\right)$.
\end{prop*}
We will first give a proof of the proposition and then continue with
the proof of the theorem. Let $\left\{ c_{n}\right\} _{n\in\mathbb{Z}}$
be in $\mathcal{N}_{s}\left(\mathbf{X},\bar{B}\right)$. Set $f\left(z\right)=\sum_{n\in\mathbb{Z}}z^{n}c_{n}$.
Then, by Lemma \ref{lem:Lemma}, the function $g:\mathbb{A^{\circ}}\rightarrow B_{0}+B_{1}$
defined by setting $g\left(s\right)=f'\left(s\right)$ and $g\left(z\right)=\frac{1}{z-s}f\left(z\right)$
for all $z\in\mathbb{A^{\circ}}\backslash\left\{ s\right\} $ has
a Laurent expansion $g\left(z\right)=\sum_{n\in\mathbb{Z}}z^{n}g_{n}$
with $\left\{ g_{n}\right\} _{n\in\mathbb{Z}}\in\mathcal{J}\left(\mathbf{X},\bar{B}\right)$.
Set $C=\max_{j=0,1}\left\Vert S\right\Vert _{\mathcal{X}_{j}\left(B_{j}\right)\rightarrow\mathcal{X}_{j}\left(B_{j}\right)}$.
Since $\mathcal{J}_{0}\left(\bar{B}\right)$ is dense in $\mathcal{J}\left(\mathbf{X},\bar{B}\right)$
we can find some $\left\{ h_{n}\right\} _{n\in\mathbb{Z}}\in\mathcal{J}_{0}\left(\bar{B}\right)$
such that \begin{equation}
\left\Vert \left\{ h_{n}\right\} _{n\in\mathbb{Z}}-\left\{ g_{n}\right\} _{n\in\mathbb{Z}}\right\Vert _{\mathcal{J}\left(\mathbf{X},\bar{B}\right)}<\frac{\varepsilon}{e\left(1+C\right)}\,.\label{eq:epsilonc}\end{equation}
For every analytic function $f:\mathbb{A}^{\circ}\rightarrow B_{0}+B_{1}$
with a Laurent expansion $f\left(z\right)=\sum_{n\in\mathbb{Z}}z^{n}b_{n}$
with $\left\{ b_{n}\right\} _{n\in\mathbb{Z}}\in\mathcal{J}\left(\mathbf{X},\bar{B}\right)$
we shall define  \[
\left\Vert f\right\Vert _{\mathcal{J},\bar{B}}:=\left\Vert \left\{ b_{n}\right\} _{n\in\mathbb{Z}}\right\Vert _{\mathcal{J}\left(\mathbf{X},\bar{B}\right)}\,.\]
(This is well-defined due to the uniqueness of the Laurent expansion
in the annulus.) Set $h\left(z\right)=\sum_{n\in\mathbb{Z}}z^{n}h_{n}$.

Note that for every element $\left\{ k_{n}\right\} _{n\in\mathbb{Z}}\in\mathcal{J}\left(\mathbf{X},\bar{B}\right)$
if $k\left(z\right)=\sum_{n\in\mathbb{Z}}z^{n}k_{n}$ and if $r\left(z\right)=z-s$
then $\left(rk\right)\left(z\right)=\sum_{n\in\mathbb{Z}}z^{n}\left(k_{n-1}-sk_{n}\right)$
and thus
\begin{eqnarray*}
\left\Vert rk\right\Vert _{\mathcal{J},\bar{B}} & = & \left\Vert \left\{ k_{n-1}-sk_{n}\right\} _{n\in\mathbb{Z}}\right\Vert _{\mathcal{J}\left(\mathbf{X},\bar{B}\right)}\\
 & \leq & \left\Vert \left\{ k_{n-1}\right\} _{n\in\mathbb{Z}}\right\Vert _{\mathcal{J}\left(\mathbf{X},\bar{B}\right)}+\left|s\right|\left\Vert \left\{ k_{n}\right\} _{n\in\mathbb{Z}}\right\Vert _{\mathcal{J}\left(\mathbf{X},\bar{B}\right)}\,.\end{eqnarray*}

Assumption (ii) of the theorem yields that \begin{eqnarray*}
\left\Vert \left\{ k_{n-1}\right\} _{n\in\mathbb{Z}}\right\Vert _{\mathcal{J}\left(\mathbf{X},\bar{B}\right)} & = & \max\left\{ \left\Vert \left\{ k_{n-1}\right\} _{n\in\mathbb{Z}}\right\Vert _{\mathcal{X}_{0}\left(B_{0}\right)},\left\Vert \left\{ e^{n}k_{n-1}\right\} _{n\in\mathbb{Z}}\right\Vert _{\mathcal{X}_{1}\left(B_{1}\right)}\right\} \\
 & \leq & C\max\left\{ \left\Vert \left\{ k_{n}\right\} _{n\in\mathbb{Z}}\right\Vert _{\mathcal{X}_{0}(B_{0})},e\left\Vert \left\{ e^{n}k_{n}\right\} _{n\in\mathbb{Z}}\right\Vert _{\mathcal{X}_{1}(B_{1})}\right\} \\
 & \le & eC\left\Vert \left\{ k_{n}\right\} _{n\in\mathbb{Z}}\right\Vert _{\mathcal{J}\left(\mathbf{X},\overline{B}\right)}\end{eqnarray*}
and thus $\left\Vert rk\right\Vert _{\mathcal{J},\bar{B}}\leq e\left(1+C\right)\left\Vert k\right\Vert _{\mathcal{J},\bar{B}}$.
The preceding calculation, along with Equation (\ref{eq:epsilonc}),
shows, in particular, that
\begin{eqnarray*}
\left\Vert \left\{ h_{n-1}-sh_{n}\right\} _{n\in\mathbb{Z}}-\left\{ c_{n}\right\} _{n\in\mathbb{Z}}\right\Vert _{\mathcal{J}\left(\mathbf{X},\bar{B}\right)} & = & \left\Vert rh-f\right\Vert _{\mathcal{J},\bar{B}}\\
 & = & \left\Vert rh-rg\right\Vert _{\mathcal{J},\bar{B}}\\
 & \le & e\left(1+C\right)\left\Vert h-g\right\Vert _{\mathcal{J},\bar{B}}<\varepsilon\,.\end{eqnarray*}
Since $\left\{ h_{n-1}-sh_{n}\right\} _{n\in\mathbb{Z}}$ is in $\mathcal{J}_{0}\left(\bar{B}\right)\cap\mathcal{N}_{s}\left(\mathbf{X},\bar{B}\right)$,
the proposition follows. 

Continuing with the proof of the theorem, we choose an element $\left\{ u_{n}\right\} _{n\in\mathbb{Z}}$
in $\mathcal{J}_{0}\left(\bar{B}\right)\cap\mathcal{N}_{s}\left(\mathbf{X},\bar{B}\right)$
such that $\left\Vert \left\{ c_{n}\right\} _{n\in\mathbb{Z}}-\left\{ u_{n}\right\} _{n\in\mathbb{Z}}\right\Vert _{\mathcal{J}\left(\mathbf{X},\bar{B}\right)}<\varepsilon/2$.

We have that $\left\{ b_{n}-u_{n}\right\} _{n\in\mathbb{Z}}\in\mathcal{J}_{0}\left(\bar{B}\right)$
and $\sum_{n\in\mathbb{Z}}s^{n}\left(b_{n}-u_{n}\right)=x$. Furthermore,\begin{eqnarray*}
 &  & \left\Vert \left\{ b_{n}-u_{n}\right\} _{n\in\mathbb{Z}}\right\Vert _{\mathcal{J}\left(\mathbf{X},\bar{B}\right)}\\
 & \le & \left\Vert \left\{ b_{n}\right\} _{n\in\mathbb{Z}}-\left\{ c_{n}\right\} _{n\in\mathbb{Z}}\right\Vert _{\mathcal{J}\left(\mathbf{X},\bar{B}\right)}+\left\Vert \left\{ c_{n}\right\} _{n\in\mathbb{Z}}-\left\{ u_{n}\right\} _{n\in\mathbb{Z}}\right\Vert _{\mathcal{J}\left(\mathbf{X},\bar{B}\right)}\\
 & < & \left\Vert x\right\Vert _{\bar{B}_{\mathbf{X},s}}+\varepsilon\,.\end{eqnarray*}
So the proof of the theorem is complete. $\qed$

\section{{\large \label{sec:Additional-remarks-regarding}Additional remarks
regarding complex interpolation}}

We begin by explicitly recalling the definition of complex interpolation
spaces on the annulus.

Let $\bar{B}=\left(B_{0},B_{1}\right)$ be a Banach pair. Let $\mathcal{F}_{\mathbb{A}}\left(\bar{B}\right)$
be the space of all continuous functions $f:\mathbb{A}\rightarrow B_{0}+B_{1}$
such that $f$ is analytic in $\mathbb{A}^{\circ}$
and for $j=0,1$ the restriction of $f$ to the circle $e^{j}\mathbb{T}$
is a continuous map of $e^{j}\mathbb{T}$ into $B_{j}$. We norm $\mathcal{F}_{\mathbb{A}}\left(\bar{B}\right)$
by $\left\Vert f\right\Vert _{\mathcal{F}_{\mathbb{A}}\left(\bar{B}\right)}=\sup\left\{ \left\Vert f\left(e^{j+it}\right)\right\Vert _{B_{j}}:t\in[0,2\pi),\, j=0,1\right\} $.
For each $\theta\in\left(0,1\right)$, let $\left[\bar{B}\right]_{\theta,\mathbb{A}}$
denote the space of all elements in $B_{0}+B_{1}$ of the form $b=f\left(e^{\theta}\right)$
for some $f\in\mathcal{F}_{\mathbb{A}}\left(\bar{B}\right)$. It is
normed by \[
\left\Vert b\right\Vert _{\left[\bar{B}\right]_{\theta,\mathbb{A}}}=\inf\left\{ \left\Vert f\right\Vert _{\mathcal{F}_{\mathbb{A}}\left(\bar{B}\right)}:f\in\mathcal{F}_{\mathbb{A}}\left(\bar{B}\right),\, b=f\left(e^{\theta}\right)\right\} \,.\]

Here, as promised above, we give a detailed proof that \[
\bar{B}_{\left\{ FC,FC\right\} ,e^{\theta}}=\left[\bar{B}\right]_{\theta,\mathbb{A}}\mbox{with equality of norms, for each }\theta\in(0,1)\,\]
as was stated in \cite{ckmr}. Some parts of the proof can also be
found on pp.~78-9 of \cite{cwikel-janson}. Related ideas appear
already in {\large \cite{calderon}}. 

Let $b\in\left[\bar{B}\right]_{\theta,\mathbb{A}}$ and let $\varepsilon$
be an arbitrary positive number. We can find some $f\in\mathcal{F}_{\mathbb{A}}\left(\bar{B}\right)$
which satisfies $b=f\left(e^{\theta}\right)$ and $\left\Vert f\right\Vert _{\mathcal{F}_{\mathbb{A}}\left(\bar{B}\right)}<\left\Vert b\right\Vert _{\left[\bar{B}\right]_{\theta,\mathbb{A}}}+\varepsilon$.
As shown on pp.~78--9 of \cite{cwikel-janson}, $f\left(z\right)=\sum_{n\in\mathbb{Z}}z^{n}\hat{f}(n)$
for every $z\in\mathbb{A^{\circ}}$, where \[
\widehat{f}(n)=\frac{1}{2\pi}\int_{0}^{2\pi}e^{-nit}f\left(e^{it}\right)dt=\frac{1}{2\pi}\int_{0}^{2\pi}e^{-n\left(1+it\right)}f\left(e^{1+it}\right)dt\,,\]
 from which we obtain that $b=\sum_{n\in\mathbb{Z}}e^{\theta n}\widehat{f}(n)$, and thus $b\in\bar{B}_{\left\{ FC,FC\right\} ,e^{\theta}}$ and
$\left[\bar{B}\right]_{\theta,\mathbb{A}}\subset\bar{B}_{\left\{ FC,FC\right\} ,e^{\theta}}\,.$
Furthermore\begin{eqnarray*}
\left\Vert b\right\Vert _{\left[\bar{B}\right]_{\theta,\mathbb{A}}}+\varepsilon & > & \sup\left\{ \left\Vert f\left(e^{j+it}\right)\right\Vert _{B_{j}}:t\in[0,2\pi),j=0,1\right\} \\
 & = & \max_{j=0,1}\left\Vert \left\{ e^{jn}\hat{f}\left(n\right)\right\} _{n\in\mathbb{Z}}\right\Vert _{FC\left(B_{j}\right)}\\
 & = & \left\Vert \left\{ \hat{f}\left(n\right)\right\} _{n\in\mathbb{Z}}\right\Vert _{\mathcal{J}\left(\left\{ FC,FC\right\} ,\bar{B}\right)}\\
 & \ge & \left\Vert b\right\Vert _{\bar{B}_{\left\{ FC,FC\right\} ,e^{\theta}}}.\end{eqnarray*}
It follows that $\left\Vert b\right\Vert _{\left[\bar{B}\right]_{\theta,\mathbb{A}}}\geq\left\Vert b\right\Vert _{\bar{B}_{\left\{ FC,FC\right\} ,e^{\theta}}}$
for all $b\in\left[\bar{B}\right]_{\theta,\mathbb{A}}\,.$

Now, for the reverse inclusion and norm inequality, let $b\in\bar{B}_{\left\{ FC,FC\right\} ,e^{\theta}}$
and let $\varepsilon$ be an arbitrary positive number. We can find
some $\left\{ b_{n}\right\} _{n\in\mathbb{Z}}\in\mathcal{J}\left(\left\{ FC,FC\right\} ,\bar{B}\right)$
such that
\begin{equation}
b=\sum_{n\in\mathbb{Z}}e^{\theta n}b_{n}\label{eq:b}\end{equation}
 and $\left\Vert \left\{ b_{n}\right\} _{n\in\mathbb{Z}}\right\Vert _{\mathcal{J}\left(\left\{ FC,FC\right\} ,\bar{B}\right)}<\left\Vert b\right\Vert _{\bar{B}_{\left\{ FC,FC\right\} ,e^{\theta}}}+\varepsilon\,.$
By the definition of $FC$, we can find for $j=0,1$ continuous functions
$f_{j}:\mathbb{T}\rightarrow B_{j}$ such that $e^{jn}b_{n}=\frac{1}{2\pi}\int_{0}^{2\pi}e^{-int}f_{j}\left(e^{it}\right)dt$
for all $n$. We shall define a sequence of functions $g_{N}:\mathbb{A}\rightarrow B_{0}\cap B_{1}$
by \[
g_{N}(z)=\sum_{n=-N}^{N}z^{n}\left(1-\frac{\left|n\right|}{N+1}\right)b_{n}\,,\]
and we shall first show that $g_{N}\left(z\right)$ converges in $B_{0}+B_{1}$
for each $z\in\mathbb{A}$. 

If $z\in\mathbb{A}^{\circ}$ then, by Remark \ref{rem:analytic},
the sequence of Laurent polynomials $S_{N}(z):=\sum_{n=-N}^{N}z^{n}b_{n}$ converges
in $B_{0}+B_{1}$, and thus, since in fact $g_{N}(z)=\frac{1}{N+1}\sum_{n=0}^{N}S_{N}(z)$,
it follows, by standard arguments, that $g_{N}(z)$ also converges
in $B_{0}+B_{1}$ to the same limit $\sum_{n\in\mathbb{Z}}z^{n}b_{n}$. 

If we apply the lemma on pp.~10--11 of \cite{katz} with the Fejér
summability kernel $K_{n}\left(t\right)=\sum_{j=-n}^{n}\left(1-\frac{\left|j\right|}{n+1}\right)e^{ijt}$
and with $\varphi\left(\tau\right)=f_{j}\left(e^{i\left(t-\tau\right)}\right)$
and $B=B_{j}$ for $j=0,1$, then we obtain that \begin{equation}
\lim_{N\rightarrow\infty}\left\Vert g_{N}\left(e^{j+it}\right)-f_{j}\left(e^{it}\right)\right\Vert _{B_{j}}=0\label{eq:limit}\end{equation}
 for each $t\in[0,2\pi)$ and thus $g_{N}\left(z\right)$ also converges
for $z\in\mathbb{A}\backslash\mathbb{A}^{\circ}$. Moreover, since
\[
\lim_{\tau\rightarrow0}\sup\left\{ \left\Vert f_{j}\left(e^{i\left(t-\tau\right)}\right)-f_{j}\left(e^{it}\right)\right\Vert _{B_{j}}:t\in[0,2\pi),\ j=0,1\right\} =0\,,\]
we can in fact obtain, by making a slight modification to the proof
of the lemma in \cite{katz} for our particular case, that
\begin{equation}
\lim_{N\rightarrow\infty}\sup\left\{ \left\Vert g_{N}\left(e^{j+it}\right)-f_{j}\left(e^{it}\right)\right\Vert _{B_{j}}:t\in[0,2\pi),\ j=0,1\right\} =0\,.\label{eq:boundary}\end{equation}
We shall denote the pointwise limit of $g_{N}$ in $B_{0}+B_{1}$
by $g$. By Equation (\ref{eq:limit}), we obtain that $g\left(e^{j+it}\right)=f_{j}\left(e^{it}\right)$
for each $t\in[0,2\pi)$ and $j=0,1$, and thus the restriction of
$g$ to the circle $e^{j}\mathbb{T}$ is a continuous map of $e^{j}\mathbb{T}$
into $B_{j}$. Since $g\left(z\right)=\sum_{n\in\mathbb{Z}}z^{n}b_{n}$
for every $z\in\mathbb{A}^{\circ}$, by Remarks \ref{rem:analytic}
and \ref{rem:admits_differentiation} and Equation (\ref{eq:b}),
$g$ is an analytic function on $\mathbb{A^{\circ}}$ which satisfies
\begin{equation}
g\left(e^{\theta}\right)=b\,.\label{eq:b2}\end{equation}
By Equation (\ref{eq:boundary}) and the maximum principle, the sequence
of continuous functions $g_{N}$
converges in $B_{0}+B_{1}$ uniformly on $\mathbb{A}$ and consequently
its limit is also a continuous $B_{0}+B_{1}$ valued function. Thus,
$g\in\mathcal{F}_{\mathbb{A}}\left(\bar{B}\right)$ and so, by Equation
(\ref{eq:b2}), $b\in\left[\bar{B}\right]_{\theta,\mathbb{A}}$ and
$\bar{B}_{\left\{ FC,FC\right\} ,e^{\theta}}\subset\left[\bar{B}\right]_{\theta,\mathbb{A}}\,.$
Furthermore, the preceding calculations show that 
\begin{eqnarray*}
\left\Vert \left\{ b_{n}\right\} _{n\in\mathbb{Z}}\right\Vert _{\mathcal{J}\left(\left\{ FC,FC\right\} ,\bar{B}\right)} & = & \max_{j=0,1}\sup_{t\in[0,2\pi)}\left\Vert f_{j}\left(e^{it}\right)\right\Vert _{B_{j}}\\
 & = & \max_{j=0,1}\sup_{t\in[0,2\pi)}\left\Vert g\left(e^{j+it}\right)\right\Vert _{B_{j}}\\
 & = & \left\Vert g\right\Vert _{\mathcal{F}_{\mathbb{A}}\left(\bar{B}\right)}\geq\left\Vert b\right\Vert _{\left[\bar{B}\right]_{\theta,\mathbb{A}}}\,.\end{eqnarray*}
Therefore, $\left\Vert b\right\Vert _{\left[\bar{B}\right]_{\theta,\mathbb{A}}}\leq\left\Vert b\right\Vert _{\bar{B}_{\left\{ FC,FC\right\} ,e^{\theta}}}$
for all $b\in\bar{B}_{\left\{ FC,FC\right\} ,e^{\theta}}\,.$ This
completes the proof.

\subsection*{Acknowledgement}

I thank Professor Michael Cwikel for several very helpful discussions
and useful comments during the preparation of this paper.


\begin{thebibliography}{HD}

\normalsize

\baselineskip=17pt

\bibitem[1]{bergh-lofstrom}J. Bergh and J. Löfström, \emph{Interpolation
spaces. An Introduction}, Springer, Berlin 1976.

\bibitem[2]{brudnyi}Y. Brudnyi, N. Krugljak, \emph{Interpolation Functors
and Interpolation Spaces}, Vol. 1, North-Holland, Amsterdam, 1991.

\bibitem[3]{calderon}A. P. Calderón, \emph{Intermediate spaces and
interpolation, the complex method}, Studia Math. 24
(1964), 113--190.

\bibitem[4]{ccrsw}R. Coifman, M. Cwikel, R. Rochberg, Y. Sagher
and G. Weiss, \emph{A theory of complex interpolation for families of Banach
spaces}, Adv. in Math. 43 (1982), 203--229.

\bibitem[5]{cwikel-complex}M. Cwikel, \emph{Complex interpolation, a
discrete definition and reiteration}, Indiana Univ. Math. J.
27 (1978), 1005--1009.

\bibitem[6]{cwikel-lecture}M. Cwikel, \emph{Lecture notes on duality
and interpolation spaces}. arXiv:0803.3558 {[}math.FA{]}.

\bibitem[7]{cwikel-janson}M.Cwikel and S. Janson, \emph{Complex interpolation
of compact operators mapping into the couple $(FL^{\infty},FL_{1}^{\infty})$},
In Contemporary Mathematics 445 (L. De Carli and M. Milman, eds.), American Mathematical Society, Providence R.I., 2007, 71--92.
Preliminary version: arXiv:math/0606551v1 {[}math.FA{]}. 

\bibitem[8]{ckmr}M. Cwikel, N. Kalton, M. Milman and R. Rochberg,
\emph{A Unified Theory of Commutator Estimates for a Class of Interpolation
Methods}, Adv. in Math. 169 (2002), 241--312.

\bibitem[9]{cwikel-peetre}M. Cwikel and J. Peetre, \emph{Abstract K
and J spaces}, J. Math. Pures Appl. 60 (1981), 1--50.

\bibitem[10]{gust-peetre}J. Gustavsson and J. Peetre, \emph{Interpolation
of Orlicz spaces}, Studia Math. 60 (1977), 33--59. 

\bibitem[11]{janson}S. Janson, \emph{Minimal and maximal methods of
interpolation}, J. Functional Analysis 44 (1981),
50--73. 

\bibitem[12]{katz}Y. Katznelson, \emph{An Introduction to Harmonic Analysis}, 
Second Corrected Edition, Dover Publications, Inc. New York (1976).

\bibitem[13]{lind-tzaf}J. Lindenstrauss and L. Tzafriri, \emph{Classical
Banach Spaces Volume I. Sequence Spaces}, Springer-Verlag, Berlin/Heidelberg/New
York (1977).

\bibitem[14]{lions-peetre}J. L. Lions and J. Peetre, \emph{Sur une
classe d'espaces d'interpolation}, Inst. Hautes Etudes Sci. Publ. Math. 19 (1964), 5--68.

\bibitem[15]{peetre-brazil}J. Peetre, \emph{A theory of interpolation
of normed spaces}, Notas de Matematica No. 39, Rio de Janeiro, 1968.
88 pp.

\bibitem[16]{peetre}J. Peetre, \emph{Sur l'utilization des suites inconditionellement
sommables dans la théorie des espaces d'interpolation}, Rend. Sem. Mat. Univ. Padova 46 (1971), 173--190.

\bibitem[17]{stafney}J. D. Stafney, \emph{The Spectrum of an Operator
on an Interpolation Space}, Trans. Amer. Math. Soc. 144 (1969), 333--349.
\end{thebibliography}
\end{document}